\documentclass[12pt]{article}
\usepackage{amsmath,amsthm}
\usepackage{amssymb,latexsym}
\usepackage{enumerate}
\usepackage[english]{babel}
\usepackage{graphicx}

\def\Z{{\mathbb Z}}

\title{
Solving binomial Thue equations 
}
\author{
Istv\'{a}n Ga\'{a}l\thanks{
        This work was partially supported by the European Union 
        and the European Social Fund through project Supercomputer, 
        the national virtual lab (grant no.: TAMOP-4.2.2.C-11/1/KONV-2012-0010)
        and also supported in part by K67580 and K75566 from the
        Hungarian National Foundation for Scientific Research,
         }\;
and L\'aszl\'o Remete
\\
University of Debrecen, Mathematical Institute \\
            H--4010 Debrecen Pf.12., Hungary \\
            e--mail: igaal@science.unideb.hu, remetel42@gmail.com \\ \\
}

\begin{document}

\maketitle
\thispagestyle{empty}

\renewcommand{\thefootnote}{}



\renewcommand{\thefootnote}{\arabic{footnote}}
\setcounter{footnote}{0}

\begin{abstract}
We consider binomial Thue equations of type $x^n-my^n=\pm 1$ in $x,y\in\Z$.
Optimizing the method of Peth\H o \cite{pet} we perform an extensive 
calculation by a high performance computer
to determine all solutions with $\max(|x|,|y|)<10^{500}$
of binomial Thue equations for $m<10^7$ for exponents 
$n=3,4,5,7,11,13,17,19,23,29$.
\end{abstract}

\section{Introduction}

The method of Peth\H o \cite{pet} (see also \cite{book}) gives a fast algorithm to calculate
"small" solutions of Thue equations. The method is based on the continued fraction 
algorithm. By "small" solutions we mean those
with absolute values less than, say $10^{500}$. Nobody belives that such 
equations have larger solutions. All our experiences show that such equations 
usually only have a few very small solutions. 

In this paper we make this algorithm more efficient in order to calculate "small" solutions 
of a special type of Thue equations, the binomial Thue equations of type
\begin{equation}
x^n-my^n=\pm 1 \;\; {\rm in}\;\; x,y\in\Z.
\label{bbb}
\end{equation}

Using sharp estimates and utilizing the specialities of these equations 
we perform an extensive calculation 
by a high performance computer to determine solutions 
with $\max(|x|,|y|)<10^{500}$ of binomial Thue equations for $1<m<10^7$
(assuming that the left hand side is irreducible)
for the exponents $n=3,4,5,7,11,13,17,19,23,29$.
Our data contains all solutions of these equations with high
probability. These results complete several results on the 
solutions of binomial Thue equations \cite{b2001}, \cite{b2010}, 
\cite{b2012}.

\section{Sharper estimates}

Let $n$ be one of $3,5,7,11,13,17,19,23,29$, in fact
our arguments are valid for any odd primes
(the case $n=4$ we shall deal with later).
Assume that for $m$ the left hand side of (\ref{bbb}) is irreducible
(we skip those $m$ for which the left hand side of (\ref{bbb}) is reducible).
Our purpose is to determine all solutions $x,y\in \Z$ of (\ref{bbb}) with
$\max(|x|,|y|)<C$. We shall perform our calculation with $C=10^{500}$.

Let $\zeta=\exp(2\pi i/n)$. Let $x,y$
be and arbitrary solution of (\ref{bbb}). Set 
\[
\beta_j=x-\zeta^{j-1} \sqrt[n]{m} \ y
\]
for $j=1,\ldots,n$, then equation (\ref{bbb}) can be written as 
\begin{equation}
\beta_1\ldots \beta_n=\pm 1.
\label{bbb2}
\end{equation}
We may assume $y\geq 0$ since on the right side we have $\pm 1$ in our equation.
Also, for $y=0$ we only have the trivial solution $x=\pm 1$.
Therefore in the following let $y\geq 1$.

For $n\neq 4$ the $\beta_2,\ldots,\beta_n$ are complex, therefore  
\[
|\beta_j|=|x-\zeta^{j-1} \sqrt[n]{m} \ y|\geq |{\rm Im} (\zeta^{j-1}) | \ \sqrt[n]{m} \ y
\]
for $j=1,\ldots, n$. This yields that 
\[
|\beta_1|=|x-\sqrt[n]{m} \ y| \leq \frac{1}{c_1 \ y^{n-1}}
\]
that is 
\begin{equation}
\left| \sqrt[n]{m}-\frac{x}{y} \right|\leq \frac{1}{c_1 \ y^n}
\label{bbb3}
\end{equation}
with
\[
c_1=(\sqrt[n]{m})^{n-1} c_2,\;\; {\rm where}\;\; c_2=\prod_{j=2}^n |{\rm Im} (\zeta^{j-1}) |.
\]
Thus we arrived at the crucial point of the method of \cite{pet}.
Obviously \mbox{$(x,y)=1$}. If the upper estimate on the right hand side of (\ref{bbb3}) satisfies
\[
\frac{1}{c_1 \ y^n}\ <\frac{1}{2y^2}
\]
that is
\begin{equation}
y> (2c_1)^{1/(n-2)},
\label{bbb4}
\end{equation}
then (\ref{bbb3}) implies
\[
\left| \sqrt[n]{m}-\frac{x}{y} \right|< \frac{1}{2 \ y^2}
\]
and appling Legenre's theorem we confer that $x/y$ is a convergent $h_j/k_j$ to $\sqrt[n]{m}$, that is
$x=h_j,y=k_j$. For the corresponding partial quotients $a_j$ it is well known that
\[
\frac{1}{(a_{j+1}+2)k_j^2}<\left| \sqrt[n]{m}-\frac{x}{y} \right|.
\]
Let $A=\max_{1\leq i\leq s}a_i$ where $k_s$ is the first denominator of a partial quotient
exceeding $C$. Combining the above estimate with (\ref{bbb3}) we get
\[
\frac{1}{(A+2)k_j^2}<\left| \sqrt[n]{m}-\frac{x}{y} \right|\leq \frac{1}{c_1 \ k_j^n}
\]
whence we obtain
\[
y=k_j< c_3= (c_1(A+2))^{n-2}.
\]
We calculate all denominators of partial quotients up to $C$, we take
their maximum, calculate the above bound and check all possible
$x=h_j,y=k_j$ for $k_j<c_3$ running again the continued fraction algorithm.

Now we return to the validity of (\ref{bbb4}).
Simple calculation shows that (\ref{bbb4}) is satisfied if
\begin{equation}
y> \left( \frac{2}{c_2}\right)^{\frac{1}{n-2}}  \ m^{-\frac{n-1}{n(n-2)}}.
\label{yy}
\end{equation}
If $m$ is large enough, then condition (\ref{yy}) is satisfied for $y\geq 1$.
Some small values of $y$ must be tested for the following values of $m$:
\[
\begin{array}{|c|c|}
\hline
n   &      \\ \hline
3		& 2\leq m\leq 4 \\ \hline
5		& 2\leq m\leq 10 \\ \hline
7		& 2\leq m\leq 29 \\ \hline
11	& 2\leq m\leq 314 \\ \hline
13	& 2\leq m\leq 1078 \\ \hline
17	& 2\leq m\leq 13489 \\ \hline
19	& 2\leq m\leq 48699 \\ \hline
23	& 2\leq m\leq 652798 \\ \hline
29	& 2\leq m\leq 7960210 \\ \hline
\end{array}
\]
For all these values of $m$ we have to test all $y$ with
\[
y \leq  \left( \frac{2}{c_2}\right)^{\frac{1}{n-2}}  \ m^{-\frac{n-1}{n(n-2)}}.
\]
Easy calculation shows that this merely yields testing $y=1$ for the values of $m$
contained in the table. This can be done very fast.

\vspace{1cm}

\noindent
{\bf Remark} \\
The case $n=4$ was considered in \cite{gr}. Remark that in that case we may assume $x>0,y>0$.
We have 
\[
|x-\sqrt[4]{m}y|\leq 1,\;\; |x\pm i\sqrt[4]{m}y|\geq \sqrt[4]{m} y,\;\;
\]
therefore
\[
|x+\sqrt[4]{m}y|\geq 2\sqrt[4]{m}y-1\geq \sqrt[4]{m}y
\]
whence
\[
\left|\sqrt[4]{m}-\frac{x}{y}\right| \leq \frac{1}{(\sqrt[4]{m})^3}\ \frac{1}{y^4}<\frac{1}{2y^2}
\]
where the last inequality is valid for all $y\geq 1$. 

\vspace{1cm}

\section{Computational aspects}

We were executing the algorithm of \cite{pet} with $C=10^{500}$ for
$n=$ 3, 4, 5, 7, 11, 13, 17, 19, 23, 29 (the case $n=4$ was dealt with in \cite{gr})
using the above sharp estimates that made the procedure for binomial Thue 
equations much more efficient. This efficient algorithm allowed us to perform
the calculations for all those $2\leq m\leq 10^7$ for which the left side of 
(\ref{bbb}) is irreducible.

The procedure was implemented in Maple \cite{maple}, we used 1200 digits accuracy. 
For each exponent this calculation involved almost $10^7$ binomial
Thue equations. The routines were running on the supercomputer (high performance computer) 
network situated in Debrecen-Budapest-P\'ecs-Szeged in Hungary under Linux. 
For each exponent the total running time was about 120-200 hours 
calculated for a single node which yields 
a few hours using parallel computing with a couple of nodes.

\section{Solutions}

In this chapter we list the results of our computation.
The triples $(m,x,y)$ in our table mean that for the $m$ there is a solution $x,y$
of equation (\ref{bbb}). The trivial solution $(m,x,y)=(m,1,0)$ is not listed. The solutions are 
displayed up to sign, that is we include only one of $(m,x,y)$ and $(m,-x,-y)$.

For each exponent $m$ we list the solutions with
\[
\max(|x|,|y|)<10^{500}
\]
of all equations with $2\leq m\leq 10^7$ for which the left hand side of 
equation (\ref{bbb}) is irreducible.

\subsection{Solutions for $n=3$}

{\linewidth=5cm

{\scriptsize

\[
\begin{array}{|c|c|c||}
\hline
 m&x& y \\ \hline 
                 2& 1& 1\\   \hline
               7& 2& 1\\   \hline
               9& 2& 1\\   \hline
              17& 18& 7\\   \hline
              19& 8& 3\\   \hline
              20& 19& 7\\   \hline
               26& 3& 1\\   \hline
              28& 3& 1\\   \hline
              37& 10& 3\\   \hline
              43& 7& 2\\   \hline
               63& 4& 1\\   \hline
              65& 4& 1\\   \hline
               91& 9& 2\\   \hline
              124& 5& 1\\   \hline
              126& 5& 1\\   \hline
             182& 17& 3\\   \hline
              215& 6& 1\\   \hline
              217& 6& 1\\   \hline
              254& 19& 3\\   \hline
              342& 7& 1\\   \hline
              344& 7& 1\\   \hline
             422& 15& 2\\   \hline
              511& 8& 1\\   \hline
              513& 8& 1\\   \hline
              614& 17& 2\\   \hline
             635& 361& 42\\   \hline
             651& 26& 3\\   \hline
              728& 9& 1\\   \hline
              730& 9& 1\\   \hline
              813& 28& 3\\   \hline
              999& 10& 1\\   \hline
             1001& 10& 1\\   \hline
             1330& 11& 1\\   \hline
             1332& 11& 1\\   \hline
             1521& 23& 2\\   \hline
             1588& 35& 3\\   \hline
             1657& 71& 6\\   \hline
             1727& 12& 1\\   \hline
             1729& 12& 1\\   \hline
             1801& 73& 6\\   \hline
             1876& 37& 3\\   \hline
             1953& 25& 2\\   \hline
             2196& 13& 1\\   \hline
             2198& 13& 1\\   \hline
             2743& 14& 1\\   \hline
             2745& 14& 1\\   \hline
\end{array}
\begin{array}{|c|c|c||}
\hline
 m&x& y \\ \hline       
             3155& 44& 3\\   \hline
             3374& 15& 1\\   \hline
             3376& 15& 1\\   \hline
             3605& 46& 3\\   \hline
             3724& 31& 2\\   \hline
             3907& 63& 4\\   \hline
             4095& 16& 1\\   \hline
             4097& 16& 1\\   \hline
             4291& 65& 4\\   \hline
             4492& 33& 2\\   \hline
             4912& 17& 1\\   \hline
             4914& 17& 1\\   \hline
            5080& 361& 21\\   \hline
             5514& 53& 3\\   \hline
             5831& 18& 1\\   \hline
             5833& 18& 1\\   \hline
             6162& 55& 3\\   \hline
             6858& 19& 1\\   \hline
             6860& 19& 1\\   \hline
             7415& 39& 2\\   \hline
             7999& 20& 1\\   \hline
             8001& 20& 1\\   \hline
             8615& 41& 2\\   \hline
             8827& 62& 3\\   \hline
             9260& 21& 1\\   \hline
             9262& 21& 1\\   \hline
             9709& 64& 3\\   \hline
             10647& 22& 1\\   \hline
            10649& 22& 1\\   \hline
             12166& 23& 1\\   \hline
            12168& 23& 1\\   \hline
            12978& 47& 2\\   \hline
            13256& 71& 3\\   \hline
            13538& 143& 6\\   \hline
             13823& 24& 1\\   \hline
            13825& 24& 1\\   \hline
            14114& 145& 6\\   \hline
             14408& 73& 3\\   \hline
             14706& 49& 2\\   \hline
            15253& 124& 5\\   \hline
             15624& 25& 1\\   \hline
            15626& 25& 1\\   \hline
            16003& 126& 5\\   \hline
            17145& 361& 14\\   \hline
             17575& 26& 1\\   \hline
            17577& 26& 1\\   \hline
          
\end{array}
\begin{array}{|c|c|c||}
\hline
 m&x& y \\ \hline               
            18745& 1036& 39\\   \hline
            18963& 80& 3\\   \hline
            19441& 242& 9\\   \hline
             19682& 27& 1\\   \hline
            19684& 27& 1\\   \hline
            19927& 244& 9\\   \hline
             20421& 82& 3\\   \hline
            20797& 55& 2\\   \hline
             21951& 28& 1\\   \hline
            21953& 28& 1\\   \hline
             23149& 57& 2\\   \hline
             24388& 29& 1\\   \hline
            24390& 29& 1\\   \hline
            26110& 89& 3\\   \hline
             26999& 30& 1\\   \hline
            27001& 30& 1\\   \hline
             27910& 91& 3\\   \hline
             29790& 31& 1\\   \hline
            29792& 31& 1\\   \hline
            31256& 63& 2\\   \hline
            32006& 127& 4\\   \hline
            32042& 667& 21\\   \hline
             32767& 32& 1\\   \hline
            32769& 32& 1\\   \hline
            33542& 129& 4\\   \hline
             34328& 65& 2\\   \hline
            34859& 98& 3\\   \hline
             35936& 33& 1\\   \hline
            35938& 33& 1\\   \hline
            37037& 100& 3\\   \hline
             39303& 34& 1\\   \hline
            39305& 34& 1\\   \hline
             42874& 35& 1\\   \hline
            42876& 35& 1\\   \hline
            44739& 71& 2\\   \hline
            45372& 107& 3\\   \hline
            46011& 215& 6\\   \hline
             46655& 36& 1\\   \hline
            46657& 36& 1\\   \hline
            47307& 217& 6\\   \hline
            47964& 109& 3\\   \hline
             48627& 73& 2\\   \hline
           48949& 4097& 112\\   \hline
             50652& 37& 1\\   \hline
            50654& 37& 1\\   \hline
             54871& 38& 1\\   \hline
\end{array}
\begin{array}{|c|c|c||}
\hline
 m&x& y \\ \hline                
            54873& 38& 1\\   \hline
            57811& 116& 3\\   \hline
             59318& 39& 1\\   \hline
            59320& 39& 1\\   \hline
            60853& 118& 3\\   \hline
            61630& 79& 2\\   \hline
             63999& 40& 1\\   \hline
            64001& 40& 1\\   \hline
             66430& 81& 2\\   \hline
             68920& 41& 1\\   \hline
            68922& 41& 1\\   \hline
            72338& 125& 3\\   \hline
             74087& 42& 1\\   \hline
            74089& 42& 1\\   \hline
            75866& 127& 3\\   \hline
             79506& 43& 1\\   \hline
            79508& 43& 1\\   \hline
            82313& 87& 2\\   \hline
             85183& 44& 1\\   \hline
            85185& 44& 1\\   \hline
          87866& 9825& 221\\   \hline
             88121& 89& 2\\   \hline
            89115& 134& 3\\   \hline
             91124& 45& 1\\   \hline
            91126& 45& 1\\   \hline
            93165& 136& 3\\   \hline
             97335& 46& 1\\   \hline
            97337& 46& 1\\   \hline
            99161& 324& 7\\   \hline
           100082& 325& 7\\   \hline
            103822& 47& 1\\   \hline
            103824& 47& 1\\   \hline
            107172& 95& 2\\   \hline
           108304& 143& 3\\   \hline
           108873& 191& 4\\   \hline
           109444& 287& 6\\   \hline
           110017& 575& 12\\   \hline
            110591& 48& 1\\   \hline
            110593& 48& 1\\   \hline
           111169& 577& 12\\   \hline
            111748& 289& 6\\   \hline
            112329& 193& 4\\   \hline
            112912& 145& 3\\   \hline
            114084& 97& 2\\   \hline
           116623& 342& 7\\   \hline
            117648& 49& 1\\   \hline
 \end{array}
\begin{array}{|c|c|c|}
\hline
 m&x& y \\ \hline               
            117650& 49& 1\\   \hline
            118681& 344& 7\\   \hline
           123506& 249& 5\\   \hline
            124999& 50& 1\\   \hline
            125001& 50& 1\\   \hline
            126506& 251& 5\\   \hline
           130067& 152& 3\\   \hline
            132650& 51& 1\\   \hline
            132652& 51& 1\\   \hline
            135269& 154& 3\\   \hline
           136591& 103& 2\\   \hline
            137160& 361& 7\\   \hline
           138303& 362& 7\\   \hline
            140607& 52& 1\\   \hline
            140609& 52& 1\\   \hline
            144703& 105& 2\\   \hline
            148876& 53& 1\\   \hline
            148878& 53& 1\\   \hline
           154566& 161& 3\\   \hline
           156494& 485& 9\\   \hline
            157463& 54& 1\\   \hline
            157465& 54& 1\\   \hline
            158438& 487& 9\\   \hline
            160398& 163& 3\\   \hline
            166374& 55& 1\\   \hline
            166376& 55& 1\\   \hline
           170954& 111& 2\\   \hline
            175615& 56& 1\\   \hline
            175617& 56& 1\\   \hline
            180362& 113& 2\\   \hline
           181963& 170& 3\\   \hline
            185192& 57& 1\\   \hline
            185194& 57& 1\\   \hline
            188461& 172& 3\\   \hline
            195111& 58& 1\\   \hline
            195113& 58& 1\\   \hline
            205378& 59& 1\\   \hline
            205380& 59& 1\\   \hline
           210645& 119& 2\\   \hline
           212420& 179& 3\\   \hline
           214205& 359& 6\\   \hline
            215999& 60& 1\\   \hline
            216001& 60& 1\\   \hline
            217805& 361& 6\\   \hline
            219620& 181& 3\\   \hline
            221445& 121& 2\\   \hline
\end{array}
\]

\[
\begin{array}{|c|c|c||}
\hline
 m&x& y \\ \hline                 
           226980& 61& 1\\   \hline
            226982& 61& 1\\   \hline
            238327& 62& 1\\   \hline
             238329& 62& 1\\   \hline
           246099& 188& 3\\   \hline
            250046& 63& 1\\   \hline
            250048& 63& 1\\   \hline
            254037& 190& 3\\   \hline
           256048& 127& 2\\   \hline
           259084& 255& 4\\   \hline
           260611& 511& 8\\   \hline
            262143& 64& 1\\   \hline
            262145& 64& 1\\   \hline
            263683& 513& 8\\   \hline
            265228& 257& 4\\   \hline
            268336& 129& 2\\   \hline
            274624& 65& 1\\   \hline
            274626& 65& 1\\   \hline
           283162& 197& 3\\   \hline
            287495& 66& 1\\   \hline
            287497& 66& 1\\   \hline
            291874& 199& 3\\   \hline
            300762& 67& 1\\   \hline
            300764& 67& 1\\   \hline
           307547& 135& 2\\   \hline
            314431& 68& 1\\   \hline
            314433& 68& 1\\   \hline
            321419& 137& 2\\   \hline
           323771& 206& 3\\   \hline
            328508& 69& 1\\   \hline
            328510& 69& 1\\   \hline
            333293& 208& 3\\   \hline
            342999& 70& 1\\   \hline
            343001& 70& 1\\   \hline
          356057& 4040& 57\\   \hline
            357910& 71& 1\\   \hline
            357912& 71& 1\\   \hline
           365526& 143& 2\\   \hline
           368088& 215& 3\\   \hline
           370662& 431& 6\\   \hline
            373247& 72& 1\\   \hline
            373249& 72& 1\\   \hline
            375846& 433& 6\\   \hline
            378456& 217& 3\\   \hline
            381078& 145& 2\\   \hline
            389016& 73& 1\\   \hline
\end{array}
\begin{array}{|c|c|c||}
\hline
 m&x& y \\ \hline   
           389018& 73& 1\\   \hline
           391592& 4097& 56\\   \hline
            405223& 74& 1\\   \hline              
            405225& 74& 1\\   \hline
           416275& 224& 3\\   \hline
           418509& 374& 5\\   \hline
          420751& 1124& 15\\   \hline
            421874& 75& 1\\   \hline
            421876& 75& 1\\   \hline
           423001& 1126& 15\\   \hline
            425259& 376& 5\\   \hline
            427525& 226& 3\\   \hline
           430369& 151& 2\\   \hline
            438975& 76& 1\\   \hline
            438977& 76& 1\\   \hline
            447697& 153& 2\\   \hline
            456532& 77& 1\\   \hline
            456534& 77& 1\\   \hline
           468494& 233& 3\\   \hline
            474551& 78& 1\\   \hline
            474553& 78& 1\\   \hline
            480662& 235& 3\\   \hline
            493038& 79& 1\\   \hline
            493040& 79& 1\\   \hline
           502460& 159& 2\\   \hline
           506115& 1036& 13\\   \hline
           507215& 319& 4\\   \hline
          507582& 1037& 13\\   \hline
            511999& 80& 1\\   \hline
            512001& 80& 1\\   \hline
            516815& 321& 4\\   \hline
            521660& 161& 2\\   \hline
           524907& 242& 3\\   \hline
           529257& 728& 9\\   \hline
            531440& 81& 1\\   \hline
            531442& 81& 1\\   \hline
            533631& 730& 9\\   \hline
            538029& 244& 3\\   \hline
            551367& 82& 1\\   \hline
            551369& 82& 1\\   \hline
            571786& 83& 1\\   \hline
            571788& 83& 1\\   \hline
           582183& 167& 2\\   \hline
           585676& 251& 3\\   \hline
           589183& 503& 6\\   \hline
            592703& 84& 1\\   \hline
\end{array}
\begin{array}{|c|c|c||}
\hline
 m&x& y \\ \hline                 
            592705& 84& 1\\   \hline
            596239& 505& 6\\   \hline
            599788& 253& 3\\   \hline
            603351& 169& 2\\   \hline
            614124& 85& 1\\   \hline
            614126& 85& 1\\   \hline
            636055& 86& 1\\   \hline
            636057& 86& 1\\   \hline
           650963& 260& 3\\   \hline
            658502& 87& 1\\   \hline
            658504& 87& 1\\   \hline
            666101& 262& 3\\   \hline
           669922& 175& 2\\   \hline
            681471& 88& 1\\   \hline
            681473& 88& 1\\   \hline
            693154& 177& 2\\   \hline
            704968& 89& 1\\   \hline
            704970& 89& 1\\   \hline
           710467& 1160& 13\\   \hline
          712306& 1161& 13\\   \hline
           720930& 269& 3\\   \hline
            728999& 90& 1\\   \hline
            729001& 90& 1\\   \hline
            737130& 271& 3\\   \hline
            753570& 91& 1\\   \hline
            753572& 91& 1\\   \hline
           766061& 183& 2\\   \hline
            778687& 92& 1\\   \hline
            778689& 92& 1\\   \hline
            791453& 185& 2\\   \hline
           795739& 278& 3\\   \hline
            804356& 93& 1\\   \hline
            804358& 93& 1\\   \hline
            813037& 280& 3\\   \hline
            830583& 94& 1\\   \hline
            830585& 94& 1\\   \hline
            857374& 95& 1\\   \hline
            857376& 95& 1\\   \hline
            865134& 667& 7\\   \hline
           869031& 668& 7\\   \hline
           870984& 191& 2\\   \hline
           875552& 287& 3\\   \hline
           877842& 383& 4\\   \hline
           880136& 575& 6\\   \hline
          882434& 1151& 12\\   \hline
            884735& 96& 1\\   \hline
\end{array}
\begin{array}{|c|c|c||}
\hline
 m&x& y \\ \hline     
           884737& 96& 1\\   \hline
           887042& 1153& 12\\   \hline
            889352& 577& 6\\   \hline            
            891666& 385& 4\\   \hline
            893984& 289& 3\\   \hline
            898632& 193& 2\\   \hline
           902629& 1353& 14\\   \hline
            912672& 97& 1\\   \hline
            912674& 97& 1\\   \hline
           937082& 685& 7\\   \hline
            941191& 98& 1\\   \hline
            941193& 98& 1\\   \hline
            945314& 687& 7\\   \hline
           960531& 296& 3\\   \hline
            970298& 99& 1\\   \hline
            970300& 99& 1\\   \hline
            980133& 298& 3\\   \hline
          980838& 1391& 14\\   \hline
           985075& 199& 2\\   \hline
           994012& 499& 5\\   \hline
           997003& 999& 10\\   \hline
            999999& 100& 1\\   \hline
           1000001& 100& 1\\   \hline
          1003003& 1001& 10\\   \hline
           1006012& 501& 5\\   \hline
           1015075& 201& 2\\   \hline
           1017241& 704& 7\\   \hline
           1021582& 705& 7\\   \hline
           1030300& 101& 1\\   \hline
           1030302& 101& 1\\   \hline
           1050838& 305& 3\\   \hline
           1061207& 102& 1\\   \hline
           1061209& 102& 1\\   \hline
           1071646& 307& 3\\   \hline
           1092726& 103& 1\\   \hline
           1092728& 103& 1\\   \hline
           1108718& 207& 2\\   \hline
           1124863& 104& 1\\   \hline
           1124865& 104& 1\\   \hline
           1141166& 209& 2\\   \hline
           1146635& 314& 3\\   \hline
           1157624& 105& 1\\   \hline
           1157626& 105& 1\\   \hline
           1168685& 316& 3\\   \hline
           1191015& 106& 1\\   \hline
           1191017& 106& 1\\   \hline
 \end{array}
\begin{array}{|c|c|c|}
\hline
 m&x& y \\ \hline               
           1225042& 107& 1\\   \hline
           1225044& 107& 1\\   \hline
           1242297& 215& 2\\   \hline
           1248084& 323& 3\\   \hline
           1253889& 647& 6\\   \hline
           1255828& 971& 9\\   \hline
          1257769& 1943& 18\\   \hline
           1259711& 108& 1\\   \hline
           1259713& 108& 1\\   \hline
          1261657& 1945& 18\\   \hline
           1263604& 973& 9\\   \hline
           1265553& 649& 6\\   \hline
           1271412& 325& 3\\   \hline
           1277289& 217& 2\\   \hline
           1295028& 109& 1\\   \hline
           1295030& 109& 1\\   \hline
          1306421& 9948& 91\\   \hline
           1330999& 110& 1\\   \hline
           1331001& 110& 1\\   \hline
           1355347& 332& 3\\   \hline
           1367630& 111& 1\\   \hline
           1367632& 111& 1\\   \hline
           1379989& 334& 3\\   \hline
           1386196& 223& 2\\   \hline
           1395541& 447& 4\\   \hline
           1404927& 112& 1\\   \hline
           1404929& 112& 1\\   \hline
           1414357& 449& 4\\   \hline
           1423828& 225& 2\\   \hline
           1442896& 113& 1\\   \hline
           1442898& 113& 1\\   \hline
           1468586& 341& 3\\   \hline
           1481543& 114& 1\\   \hline
           1481545& 114& 1\\   \hline
           1494578& 343& 3\\   \hline
           1520874& 115& 1\\   \hline
           1520876& 115& 1\\   \hline
          1530341& 2420& 21\\   \hline
           1540799& 231& 2\\   \hline
           1560895& 116& 1\\   \hline
           1560897& 116& 1\\   \hline
           1581167& 233& 2\\   \hline
           1587963& 350& 3\\   \hline
           1601612& 117& 1\\   \hline
           1601614& 117& 1\\   \hline
           1615341& 352& 3\\   \hline
\end{array}
\]

\[
\begin{array}{|c|c|c||}
\hline
 m&x& y \\ \hline                
           1643031& 118& 1\\   \hline
           1643033& 118& 1\\   \hline
           1685158& 119& 1\\   \hline
           1685160& 119& 1\\   \hline
          1695638& 7751& 65\\   \hline
           1706490& 239& 2\\   \hline
           1713640& 359& 3\\   \hline
           1720810& 719& 6\\   \hline
           1727999& 120& 1\\   \hline
           1728001& 120& 1\\   \hline
           1735210& 721& 6\\   \hline
           1742440& 361& 3\\   \hline
           1749690& 241& 2\\   \hline
          1767571& 1330& 11\\   \hline
           1771560& 121& 1\\   \hline
           1771562& 121& 1\\   \hline
          1775557& 1332& 11\\   \hline
           1815847& 122& 1\\   \hline
           1815849& 122& 1\\   \hline
           1845779& 368& 3\\   \hline
           1860866& 123& 1\\   \hline
           1860868& 123& 1\\   \hline
           1876037& 370& 3\\   \hline
           1883653& 247& 2\\   \hline
           1906623& 124& 1\\   \hline
           1906625& 124& 1\\   \hline
          1922636& 3233& 26\\   \hline
           1929781& 249& 2\\   \hline
           1943765& 624& 5\\   \hline
           1953124& 125& 1\\   \hline
           1953126& 125& 1\\   \hline
           1962515& 626& 5\\   \hline
           1984542& 377& 3\\   \hline
           2000375& 126& 1\\   \hline
           2000377& 126& 1\\   \hline
           2016294& 379& 3\\   \hline
           2048382& 127& 1\\   \hline
           2048384& 127& 1\\   \hline
           2072672& 255& 2\\   \hline
           2084888& 511& 4\\   \hline
          2091014& 1023& 8\\   \hline
           2097151& 128& 1\\   \hline
           2097153& 128& 1\\   \hline
           2103302& 1025& 8\\   \hline
           2109464& 513& 4\\   \hline
           2121824& 257& 2\\   \hline
\end{array}
\begin{array}{|c|c|c||}
\hline
 m&x& y \\ \hline       
        2130091& 386& 3\\   \hline
           2146688& 129& 1\\   \hline         
           2146690& 129& 1\\   \hline
           2163373& 388& 3\\   \hline
          2187357& 2726& 21\\   \hline
           2196999& 130& 1\\   \hline
           2197001& 130& 1\\   \hline
           2248090& 131& 1\\   \hline
           2248092& 131& 1\\   \hline
           2273931& 263& 2\\   \hline
           2282588& 395& 3\\   \hline
           2291267& 791& 6\\   \hline
           2299967& 132& 1\\   \hline
           2299969& 132& 1\\   \hline
           2308691& 793& 6\\   \hline
           2317436& 397& 3\\   \hline
           2326203& 265& 2\\   \hline
           2352636& 133& 1\\   \hline
           2352638& 133& 1\\   \hline
           2406103& 134& 1\\   \hline
           2406105& 134& 1\\   \hline
           2442195& 404& 3\\   \hline
          2454305& 1214& 9\\   \hline
           2460374& 135& 1\\   \hline
           2460376& 135& 1\\   \hline
           2466455& 1216& 9\\   \hline
           2478645& 406& 3\\   \hline
           2487814& 271& 2\\   \hline
           2515455& 136& 1\\   \hline
           2515457& 136& 1\\   \hline
           2543302& 273& 2\\   \hline
           2571352& 137& 1\\   \hline
           2571354& 137& 1\\   \hline
           2609074& 413& 3\\   \hline
           2628071& 138& 1\\   \hline
           2628073& 138& 1\\   \hline
           2647162& 415& 3\\   \hline
           2685618& 139& 1\\   \hline
           2685620& 139& 1\\   \hline
           2714705& 279& 2\\   \hline
           2743999& 140& 1\\   \hline
           2744001& 140& 1\\   \hline
           2773505& 281& 2\\   \hline
           2783387& 422& 3\\   \hline
           2803220& 141& 1\\   \hline
           2803222& 141& 1\\   \hline
\end{array}
\begin{array}{|c|c|c||}
\hline
 m&x& y \\ \hline               
          2823149& 424& 3\\   \hline
           2863287& 142& 1\\   \hline 
           2863289& 142& 1\\   \hline
          2888177& 5554& 39\\   \hline
           2924206& 143& 1\\   \hline
           2924208& 143& 1\\   \hline
           2954988& 287& 2\\   \hline
           2965296& 431& 3\\   \hline
           2970459& 575& 4\\   \hline
           2975628& 863& 6\\   \hline
          2980803& 1727& 12\\   \hline
           2985983& 144& 1\\   \hline
           2985985& 144& 1\\   \hline
          2991171& 1729& 12\\   \hline
           2996364& 865& 6\\   \hline
           3001563& 577& 4\\   \hline
           3003793& 1010& 7\\   \hline
           3006768& 433& 3\\   \hline
          3012724& 1011& 7\\   \hline
           3017196& 289& 2\\   \hline
           3048624& 145& 1\\   \hline
           3048626& 145& 1\\   \hline
           3112135& 146& 1\\   \hline
           3112137& 146& 1\\   \hline
          3132736& 4097& 28\\   \hline
          3141461& 5126& 35\\   \hline
           3154963& 440& 3\\   \hline
          3167271& 1028& 7\\   \hline
          3173437& 3086& 21\\   \hline
           3176522& 147& 1\\   \hline
           3176524& 147& 1\\   \hline
          3179611& 3088& 21\\   \hline
           3185793& 1030& 7\\   \hline
           3198181& 442& 3\\   \hline
           3209047& 295& 2\\   \hline
           3241791& 148& 1\\   \hline
           3241793& 148& 1\\   \hline
          3266062& 2819& 19\\   \hline
          3269539& 2820& 19\\   \hline
           3274759& 297& 2\\   \hline
           3307948& 149& 1\\   \hline
           3307950& 149& 1\\   \hline
           3346154& 1047& 7\\   \hline
           3352550& 449& 3\\   \hline
           3355751& 1048& 7\\   \hline
           3361518& 749& 5\\   \hline
\end{array}
\begin{array}{|c|c|c||}
\hline
 m&x& y \\ \hline                
           3370502& 2249& 15\\   \hline
           3374999& 150& 1\\   \hline
           3375001& 150& 1\\   \hline
          3379502& 2251& 15\\   \hline
           3388518& 751& 5\\   \hline
           3397550& 451& 3\\   \hline
           3442950& 151& 1\\   \hline
           3442952& 151& 1\\   \hline
           3477266& 303& 2\\   \hline
           3511807& 152& 1\\   \hline
           3511809& 152& 1\\   \hline
           3546578& 305& 2\\   \hline
           3558219& 458& 3\\   \hline
           3581576& 153& 1\\   \hline
           3581578& 153& 1\\   \hline
           3605037& 460& 3\\   \hline
           3652263& 154& 1\\   \hline
           3652265& 154& 1\\   \hline
           3723874& 155& 1\\   \hline
           3723876& 155& 1\\   \hline
           3760029& 311& 2\\   \hline
           3772132& 467& 3\\   \hline
           3784261& 935& 6\\   \hline
           3796415& 156& 1\\   \hline
           3796417& 156& 1\\   \hline
           3808597& 937& 6\\   \hline
           3820804& 469& 3\\   \hline
           3833037& 313& 2\\   \hline
           3869892& 157& 1\\   \hline
           3869894& 157& 1\\   \hline
           3944311& 158& 1\\   \hline
           3944313& 158& 1\\   \hline
           3994451& 476& 3\\   \hline
           4019678& 159& 1\\   \hline
           4019680& 159& 1\\   \hline
           4045013& 478& 3\\   \hline
           4057720& 319& 2\\   \hline
           4076830& 639& 4\\   \hline
           4095999& 160& 1\\   \hline
           4096001& 160& 1\\   \hline
           4115230& 641& 4\\   \hline
           4134520& 321& 2\\   \hline
           4173280& 161& 1\\   \hline
           4173282& 161& 1\\   \hline
           4225338& 485& 3\\   \hline
          4242786& 1457& 9\\   \hline
 \end{array}
\begin{array}{|c|c|c|}
\hline
 m&x& y \\ \hline    
            4251527& 162& 1\\   \hline   
            4251529& 162& 1\\   \hline        
           4260282& 1459& 9\\   \hline
           4277826& 487& 3\\   \hline
          4321265& 6841& 42\\   \hline
           4330746& 163& 1\\   \hline
           4330748& 163& 1\\   \hline
           4370723& 327& 2\\   \hline
           4410943& 164& 1\\   \hline
           4410945& 164& 1\\   \hline
          4426331& 3448& 21\\   \hline
           4451411& 329& 2\\   \hline
           4464955& 494& 3\\   \hline
           4492124& 165& 1\\   \hline
           4492126& 165& 1\\   \hline
           4519405& 496& 3\\   \hline
           4574295& 166& 1\\   \hline
           4574297& 166& 1\\   \hline
           4657462& 167& 1\\   \hline
           4657464& 167& 1\\   \hline
           4699422& 335& 2\\   \hline
           4713464& 503& 3\\   \hline
          4727534& 1007& 6\\   \hline
           4741631& 168& 1\\   \hline
           4741633& 168& 1\\   \hline
           4755758& 1009& 6\\   \hline
           4769912& 505& 3\\   \hline
           4784094& 337& 2\\   \hline
          4820221& 2196& 13\\   \hline
           4826808& 169& 1\\   \hline
           4826810& 169& 1\\   \hline
          4833403& 2198& 13\\   \hline
           4912999& 170& 1\\   \hline
           4913001& 170& 1\\   \hline
          4931602& 2383& 14\\   \hline
           4971027& 512& 3\\   \hline
           5000210& 171& 1\\   \hline
           5000212& 171& 1\\   \hline
           5029509& 514& 3\\   \hline
           5044201& 343& 2\\   \hline
           5088447& 172& 1\\   \hline
           5088449& 172& 1\\   \hline
           5132953& 345& 2\\   \hline
           5177716& 173& 1\\   \hline
           5177718& 173& 1\\   \hline
           5237806& 521& 3\\   \hline
\end{array}
\]

\[
\begin{array}{|c|c|c||}
\hline
 m&x& y \\ \hline       
           5268023& 174& 1\\   \hline
           5268025& 174& 1\\   \hline         
           5298358& 523& 3\\   \hline
           5341021& 874& 5\\   \hline
           5359374& 175& 1\\   \hline
           5359376& 175& 1\\   \hline
           5377771& 876& 5\\   \hline
           5405444& 351& 2\\   \hline
           5428577& 703& 4\\   \hline
           5451775& 176& 1\\   \hline
           5451777& 176& 1\\   \hline
           5475041& 705& 4\\   \hline
           5498372& 353& 2\\   \hline
           5513963& 530& 3\\   \hline
           5545232& 177& 1\\   \hline
           5545234& 177& 1\\   \hline
           5576621& 532& 3\\   \hline
           5639751& 178& 1\\   \hline
           5639753& 178& 1\\   \hline
          5712483& 3754& 21\\   \hline
           5735338& 179& 1\\   \hline
           5735340& 179& 1\\   \hline
           5783535& 359& 2\\   \hline
           5799660& 539& 3\\   \hline
          5815815& 1079& 6\\   \hline
           5831999& 180& 1\\   \hline
           5832001& 180& 1\\   \hline
           5848215& 1081& 6\\   \hline
           5864460& 541& 3\\   \hline
           5880735& 361& 2\\   \hline
           5929740& 181& 1\\   \hline
           5929742& 181& 1\\   \hline
           6028567& 182& 1\\   \hline
           6028569& 182& 1\\   \hline
           6095059& 548& 3\\   \hline
           6128486& 183& 1\\   \hline
           6128488& 183& 1\\   \hline
           6162037& 550& 3\\   \hline
           6178858& 367& 2\\   \hline
           6229503& 184& 1\\   \hline
           6229505& 184& 1\\   \hline
           6280426& 369& 2\\   \hline
           6331624& 185& 1\\   \hline
           6331626& 185& 1\\   \hline
           6400322& 557& 3\\   \hline
          6402292& 6499& 35\\   \hline
\end{array}
\begin{array}{|c|c|c||}
\hline
 m&x& y \\ \hline       
           6434855& 186& 1\\   \hline
           6434857& 186& 1\\   \hline         
           6469514& 559& 3\\   \hline
           6539202& 187& 1\\   \hline
           6539204& 187& 1\\   \hline
           6591797& 375& 2\\   \hline
           6644671& 188& 1\\   \hline
           6644673& 188& 1\\   \hline
           6697829& 377& 2\\   \hline
           6715611& 566& 3\\   \hline
          6739369& 1700& 9\\   \hline
           6751268& 189& 1\\   \hline
           6751270& 189& 1\\   \hline
           6763183& 1702& 9\\   \hline
           6787053& 568& 3\\   \hline
           6858999& 190& 1\\   \hline
           6859001& 190& 1\\   \hline
           6967870& 191& 1\\   \hline
           6967872& 191& 1\\   \hline
           7022736& 383& 2\\   \hline
           7041088& 575& 3\\   \hline
           7050276& 767& 4\\   \hline
          7059472& 1151& 6\\   \hline
          7064073& 1535& 8\\   \hline
          7068676& 2303& 12\\   \hline
          7073281& 4607& 24\\   \hline
           7077887& 192& 1\\   \hline
           7077889& 192& 1\\   \hline
          7082497& 4609& 24\\   \hline
          7087108& 2305& 12\\   \hline
           7091721& 1537& 8\\   \hline
           7096336& 1153& 6\\   \hline
           7105572& 769& 4\\   \hline
           7114816& 577& 3\\   \hline
           7133328& 385& 2\\   \hline
           7189056& 193& 1\\   \hline
           7189058& 193& 1\\   \hline
           7221032& 1353& 7\\   \hline
          7237055& 1354& 7\\   \hline
           7301383& 194& 1\\   \hline
           7301385& 194& 1\\   \hline
           7376915& 584& 3\\   \hline
           7414874& 195& 1\\   \hline
           7414876& 195& 1\\   \hline
           7452965& 586& 3\\   \hline
           7472059& 391& 2\\   \hline
 \end{array}
\begin{array}{|c|c|c||}
\hline
 m&x& y \\ \hline              
            7513084& 1371& 7\\   \hline
          7521307& 2743& 14\\   \hline
           7529535& 196& 1\\   \hline
           7529537& 196& 1\\   \hline
          7537771& 2745& 14\\   \hline
           7546012& 1373& 7\\   \hline
           7587307& 393& 2\\   \hline
           7645372& 197& 1\\   \hline
           7645374& 197& 1\\   \hline
           7723254& 593& 3\\   \hline
           7762391& 198& 1\\   \hline
           7762393& 198& 1\\   \hline
           7801662& 595& 3\\   \hline
           7829793& 1390& 7\\   \hline
          7846704& 1391& 7\\   \hline
           7880598& 199& 1\\   \hline
           7880600& 199& 1\\   \hline
         7926126& 93311& 468\\   \hline
           7940150& 399& 2\\   \hline
           7976024& 999& 5\\   \hline
          7988006& 1999& 10\\   \hline
           7999999& 200& 1\\   \hline
           8000001& 200& 1\\   \hline
          8012006& 2001& 10\\   \hline
           8024024& 1001& 5\\   \hline
           8060150& 401& 2\\   \hline
           8080267& 602& 3\\   \hline
           8120600& 201& 1\\   \hline
           8120602& 201& 1\\   \hline
           8161069& 604& 3\\   \hline
           8242407& 202& 1\\   \hline
           8242409& 202& 1\\   \hline
          8341522& 6287& 31\\   \hline
          8345503& 6288& 31\\   \hline
           8365426& 203& 1\\   \hline
           8365428& 203& 1\\   \hline
           8427393& 407& 2\\   \hline
           8448116& 611& 3\\   \hline
          8468873& 1223& 6\\   \hline
           8489663& 204& 1\\   \hline
           8489665& 204& 1\\   \hline
           8510489& 1225& 6\\   \hline
           8531348& 613& 3\\   \hline
           8552241& 409& 2\\   \hline
           8615124& 205& 1\\   \hline
           8615126& 205& 1\\   \hline
\end{array}
\begin{array}{|c|c|c|}
\hline
 m&x& y \\ \hline                
            8741815& 206& 1\\   \hline
           8741817& 206& 1\\   \hline
           8826963& 620& 3\\   \hline
           8869742& 207& 1\\   \hline
           8869744& 207& 1\\   \hline
           8912661& 622& 3\\   \hline
           8934172& 415& 2\\   \hline
           8966503& 831& 4\\   \hline
           8998911& 208& 1\\   \hline
           8998913& 208& 1\\   \hline
           9031399& 833& 4\\   \hline
           9063964& 417& 2\\   \hline
          9114217& 5431& 26\\   \hline
           9129328& 209& 1\\   \hline
           9129330& 209& 1\\   \hline
           9216970& 629& 3\\   \hline
           9260999& 210& 1\\   \hline
           9261001& 210& 1\\   \hline
           9305170& 631& 3\\   \hline
           9393930& 211& 1\\   \hline
           9393932& 211& 1\\   \hline
           9460871& 423& 2\\   \hline
           9528127& 212& 1\\   \hline
           9528129& 212& 1\\   \hline
           9595703& 425& 2\\   \hline
         9603769& 17855& 84\\   \hline
          9606402& 4039& 19\\   \hline
          9613539& 4040& 19\\   \hline
           9618299& 638& 3\\   \hline
           9663596& 213& 1\\   \hline
           9663598& 213& 1\\   \hline
           9709037& 640& 3\\   \hline
           9800343& 214& 1\\   \hline
           9800345& 214& 1\\   \hline
           9938374& 215& 1\\   \hline
           9938376& 215& 1\\   \hline
\end{array}
\]
}

}

\subsection{Solutions for $n=4$}

{\scriptsize
\[
\begin{array}{|c|c|c||}
\hline
 m&x& y \\ \hline 
        2  &   1  &   1  \\  \hline    
        5  &   3  &   2  \\  \hline        
        15  &   2  &   1  \\  \hline        
        17  &   2  &   1  \\  \hline        
       39  &   5  &   2  \\  \hline          
        80  &   3  &   1  \\  \hline          
        82  &   3  &   1  \\  \hline           
      150  &   7  &   2  \\  \hline            
       255  &   4  &   1  \\  \hline           
       257  &   4  &   1  \\  \hline             
      410  &   9  &   2  \\  \hline               
       624  &   5  &   1  \\  \hline              
       626  &   5  &   1  \\  \hline              
     915  &   11  &   2  \\  \hline              
      1295  &   6  &   1  \\  \hline              
      1297  &   6  &   1  \\  \hline               
     1785  &   13  &   2  \\  \hline              
      2400  &   7  &   1  \\  \hline               
      2402  &   7  &   1  \\  \hline                
     3164  &   15  &   2  \\  \hline                 
      4095  &   8  &   1  \\  \hline                  
      4097  &   8  &   1  \\  \hline                   
     5220  &   17  &   2  \\  \hline                   
      6560  &   9  &   1  \\  \hline                     
      6562  &   9  &   1  \\  \hline                     
     7140  &   239  &   26  \\  \hline                 
     8145  &   19  &   2  \\  \hline                       
     9999  &   10  &   1  \\  \hline                      
    10001  &   10  &   1  \\  \hline                   
    12155  &   21  &   2  \\  \hline                  
    14640  &   11  &   1  \\  \hline                
    14642  &   11  &   1  \\  \hline                 
    17490  &   23  &   2  \\  \hline                   
    20735  &   12  &   1  \\  \hline                 
    20737  &   12  &   1  \\  \hline                 
    24414  &   25  &   2  \\  \hline                    
    28560  &   13  &   1  \\  \hline                  
    28562  &   13  &   1  \\  \hline                  
    33215  &   27  &   2  \\  \hline                 
    38415  &   14  &   1  \\  \hline                   
    38417  &   14  &   1  \\  \hline                    
    44205  &   29  &   2  \\  \hline                   
    50624  &   15  &   1  \\  \hline                   
    50626  &   15  &   1  \\  \hline                    
    57720  &   31  &   2  \\  \hline                     
   61535  &   63  &   4  \\  \hline                           
\end{array}  
\begin{array}{|c|c|c||}
\hline
 m&x& y \\ \hline 
    65535  &   16  &   1  \\  \hline                     
    65537  &   16  &   1  \\  \hline                     
   69729  &   65  &   4  \\  \hline                    
    74120  &   33  &   2  \\  \hline     
    83520  &   17  &   1  \\  \hline                     
    83522  &   17  &   1  \\  \hline                    
    93789  &   35  &   2  \\  \hline                   
   104975  &   18  &   1  \\  \hline                   
    104977  &   18  &   1  \\  \hline                   
114240  &   239  &   13  \\  \hline               
   117135  &   37  &   2  \\  \hline                   
    130320  &   19  &   1  \\  \hline                 
    130322  &   19  &   1  \\  \hline                
   144590  &   39  &   2  \\  \hline                  
    159999  &   20  &   1  \\  \hline              
    160001  &   20  &   1  \\  \hline                    
   176610  &   41  &   2  \\  \hline                     
    194480  &   21  &   1  \\  \hline                  
    194482  &   21  &   1  \\  \hline                  
   213675  &   43  &   2  \\  \hline                   
    234255  &   22  &   1  \\  \hline                   
    234257  &   22  &   1  \\  \hline                    
   256289  &   45  &   2  \\  \hline                  
    279840  &   23  &   1  \\  \hline                     
    279842  &   23  &   1  \\  \hline                    
   304980  &   47  &   2  \\  \hline                    
    331775  &   24  &   1  \\  \hline                      
    331777  &   24  &   1  \\  \hline                 
   360300  &   49  &   2  \\  \hline                    
    390624  &   25  &   1  \\  \hline                  
    390626  &   25  &   1  \\  \hline                    
  422825  &   51  &   2  \\  \hline                      
    456975  &   26  &   1  \\  \hline                    
    456977  &   26  &   1  \\  \hline                    
   493155  &   53  &   2  \\  \hline                  
   505679  &   80  &   3  \\  \hline                    
518440  &   161  &   6  \\  \hline                   
    531440  &   27  &   1  \\  \hline                   
    531442  &   27  &   1  \\  \hline                     
544685  &   163  &   6  \\  \hline                       
   558175  &   82  &   3  \\  \hline                     
   571914  &   55  &   2  \\  \hline                     
    614655  &   28  &   1  \\  \hline                  
    614657  &   28  &   1  \\  \hline                   
   659750  &   57  &   2  \\  \hline                   
   707280  &   29  &   1  \\  \hline                      
\end{array}  
\begin{array}{|c|c|c||}
\hline
 m&x& y \\ \hline 
     707282  &   29  &   1  \\  \hline                     
   757335  &   59  &   2  \\  \hline                       
    809999  &   30  &   1  \\  \hline                     
    810001  &   30  &   1  \\  \hline                      
   865365  &   61  &   2  \\  \hline                     
    923520  &   31  &   1  \\  \hline                     
    923522  &   31  &   1  \\  \hline    
     984560  &   63  &   2  \\  \hline                       
   1016190  &   127  &   4  \\  \hline                    
   1048575  &   32  &   1  \\  \hline                    
   1048577  &   32  &   1  \\  \hline                   
1081730  &   129  &   4  \\  \hline                    
1115664  &   65  &   2  \\  \hline                       
   1185920  &   33  &   1  \\  \hline                       
   1185922  &   33  &   1  \\  \hline                   
  1259445  &   67  &   2  \\  \hline                       
   1336335  &   34  &   1  \\  \hline                        
   1336337  &   34  &   1  \\  \hline                       
  1416695  &   69  &   2  \\  \hline                   
   1500624  &   35  &   1  \\  \hline                      
   1500626  &   35  &   1  \\  \hline                        
  1588230  &   71  &   2  \\  \hline                         
  1679615  &   36  &   1  \\  \hline                        
   1679617  &   36  &   1  \\  \hline                  
1755519  &   182  &   5  \\  \hline                       
  1774890  &   73  &   2  \\  \hline                      
   1874160  &   37  &   1  \\  \hline                       
   1874162  &   37  &   1  \\  \hline                        
  1977539  &   75  &   2  \\  \hline                         
  2085135  &   38  &   1  \\  \hline                     
   2085137  &   38  &   1  \\  \hline                       
2197065  &   77  &   2  \\  \hline                           
  2313440  &   39  &   1  \\  \hline                       
   2313442  &   39  &   1  \\  \hline                      
  2434380  &   79  &   2  \\  \hline                     
   2559999  &   40  &   1  \\  \hline                       
   2560001  &   40  &   1  \\  \hline                     
2690420  &   81  &   2  \\  \hline                      
  2825760  &   41  &   1  \\  \hline                    
   2825762  &   41  &   1  \\  \hline                 
2966145  &   83  &   2  \\  \hline                     
  3111695  &   42  &   1  \\  \hline                          
   3111697  &   42  &   1  \\  \hline                          
  3262539  &   85  &   2  \\  \hline                        
   3418800  &   43  &   1  \\  \hline                        
   3418802  &   43  &   1  \\  \hline                                    
\end{array}  
\begin{array}{|c|c|c|}
\hline
 m&x& y \\ \hline 
  3580610  &   87  &   2  \\  \hline      
  3748095  &   44  &   1  \\  \hline                     
   3748097  &   44  &   1  \\  \hline                    
3851367  &   443  &   10  \\  \hline                    
3921390  &   89  &   2  \\  \hline                         
   4100624  &   45  &   1  \\  \hline                       
   4100626  &   45  &   1  \\  \hline                   
4285935  &   91  &   2  \\  \hline                        
  4477455  &   46  &   1  \\  \hline                      
   4477457  &   46  &   1  \\  \hline                    
  4675325  &   93  &   2  \\  \hline   
     4879680  &   47  &   1  \\  \hline                     
   4879682  &   47  &   1  \\  \hline                         
  5090664  &   95  &   2  \\  \hline                      
5198685  &   191  &   4  \\  \hline                   
   5308415  &   48  &   1  \\  \hline                   
   5308417  &   48  &   1  \\  \hline                     
5419875  &   193  &   4  \\  \hline                     
  5533080  &   97  &   2  \\  \hline                       
   5764800  &   49  &   1  \\  \hline                     
   5764802  &   49  &   1  \\  \hline                    
6003725  &   99  &   2  \\  \hline                       
   6249999  &   50  &   1  \\  \hline                  
   6250001  &   50  &   1  \\  \hline                   
6503775  &   101  &   2  \\  \hline                    
  6765200  &   51  &   1  \\  \hline                     
   6765202  &   51  &   1  \\  \hline                    
7034430  &   103  &   2  \\  \hline                      
    7311615  &   52  &   1  \\  \hline                   
  7311617  &   52  &   1  \\  \hline                      
7596914  &   105  &   2  \\  \hline                   
   7890480  &   53  &   1  \\  \hline                     
   7890482  &   53  &   1  \\  \hline                   
8192475  &   107  &   2  \\  \hline                   
8295040  &   161  &   3  \\  \hline                    
8398565  &   323  &   6  \\  \hline                 
   8503055  &   54  &   1  \\  \hline                    
   8503057  &   54  &   1  \\  \hline                    
8608519  &   325  &   6  \\  \hline                     
8714960  &   163  &   3  \\  \hline                      
8822385  &   109  &   2  \\  \hline                       
   9150624  &   55  &   1  \\  \hline                     
   9150626  &   55  &   1  \\  \hline                     
9487940  &   111  &   2  \\  \hline                    
  9834495  &   56  &   1  \\  \hline                       
   9834497  &   56  &   1  \\  \hline                      
\end{array}
\]
}

\subsection{Solutions for $n=5$}

{\scriptsize
\[
\begin{array}{|c|c|c||}
\hline
 m&x& y \\ \hline 
                2& 1& 1\\ \hline
                31& 2& 1\\ \hline
               33& 2& 1\\ \hline
               242& 3& 1\\ \hline
               244& 3& 1\\ \hline
               1023& 4& 1\\ \hline
              1025& 4& 1\\ \hline
               3124& 5& 1\\ \hline
              3126& 5& 1\\ \hline
               7775& 6& 1\\ \hline
              7777& 6& 1\\ \hline
              16806& 7& 1\\ \hline 
              16808& 7& 1\\ \hline
\end{array}   
\begin{array}{|c|c|c||}
\hline
 m&x& y \\ \hline 
              32767& 8& 1\\ \hline
              32769& 8& 1\\ \hline
              59048& 9& 1\\ \hline
              59050& 9& 1\\ \hline
              99999& 10& 1\\ \hline
             100001& 10& 1\\ \hline
             161050& 11& 1\\ \hline
             161052& 11& 1\\ \hline
             248831& 12& 1\\ \hline
             248833& 12& 1\\ \hline
             371292& 13& 1\\ \hline
             371294& 13& 1\\ \hline
             537823& 14& 1\\ \hline
\end{array}
\begin{array}{|c|c|c||}
\hline
 m&x& y \\ \hline            
             537825& 14& 1\\ \hline
             759374& 15& 1\\ \hline
             759376& 15& 1\\ \hline
             894661& 31& 2\\ \hline
             1048575& 16& 1\\ \hline
            1048577& 16& 1\\ \hline
             1222981& 33& 2\\ \hline
             1419856& 17& 1\\ \hline
            1419858& 17& 1\\ \hline
             1889567& 18& 1\\ \hline
            1889569& 18& 1\\ \hline
             2476098& 19& 1\\ \hline
            2476100& 19& 1\\ \hline
\end{array}
\begin{array}{|c|c|c|}
\hline
 m&x& y \\ \hline 
             3199999& 20& 1\\ \hline
            3200001& 20& 1\\ \hline
             4084100& 21& 1\\ \hline
            4084102& 21& 1\\ \hline
             5153631& 22& 1\\ \hline
            5153633& 22& 1\\ \hline
             6436342& 23& 1\\ \hline
            6436344& 23& 1\\ \hline
             7962623& 24& 1\\ \hline
            7962625& 24& 1\\ \hline
            9765624& 25& 1\\ \hline
            9765626& 25& 1\\ \hline
\end{array}
\]
}

\subsection{Solutions for $n=7$}

{\scriptsize
\[
\begin{array}{|c|c|c|}
\hline
 m&x& y \\ \hline 
                2& 1& 1\\ \hline
               127& 2& 1\\ \hline
               129& 2& 1\\ \hline
               2186& 3& 1\\ \hline
              2188& 3& 1\\ \hline
              16383& 4& 1\\ \hline
              16385& 4& 1\\ \hline
              78124& 5& 1\\ \hline
              78126& 5& 1\\ \hline
              279935& 6& 1\\ \hline
             279937& 6& 1\\ \hline
              823542& 7& 1\\ \hline
             823544& 7& 1\\ \hline
             2097151& 8& 1\\ \hline
             2097153& 8& 1\\ \hline
             4782968& 9& 1\\ \hline
             4782970& 9& 1\\ \hline
            9999999& 10& 1\\ \hline
\end{array}
\]
}

\subsection{Solutions for $n=11$}

{\scriptsize
\[
\begin{array}{|c|c|c|}
\hline
 m&x& y \\ \hline   
                2& 1& 1\\ \hline
               2047& 2& 1\\ \hline
              2049& 2& 1\\ \hline
              177146& 3& 1\\ \hline
             177148& 3& 1\\ \hline
             4194303& 4& 1\\ \hline
             4194305& 4& 1\\ \hline
\end{array}
\]
}

\subsection{Solutions for $n=13$}

{\scriptsize
\[
\begin{array}{|c|c|c|}
\hline
 m&x& y \\ \hline 
                2& 1& 1\\ \hline
               8191& 2& 1\\ \hline
              8193& 2& 1\\ \hline
             1594322& 3& 1\\ \hline
             1594324& 3& 1\\ \hline
\end{array}
\]
}

\subsection{Solutions for $n=17$}

{\scriptsize
\[
\begin{array}{|c|c|c|}
\hline
 m&x& y \\ \hline 
                2& 1& 1\\ \hline
              131071& 2& 1\\ \hline
             131073& 2& 1\\ \hline
\end{array}
\]
}

\subsection{Solutions for $n=19$}

{\scriptsize
\[
\begin{array}{|c|c|c|}
\hline
 m&x& y \\ \hline 
                2& 1& 1\\ \hline
              524287& 2& 1\\ \hline
             524289& 2& 1\\ \hline
\end{array}
\]
}

\subsection{Solutions for $n=23$}

{\scriptsize
\[
\begin{array}{|c|c|c|}
\hline
 m&x& y \\ \hline 
                2& 1& 1\\ \hline
             8388607& 2& 1\\ \hline
             8388609& 2& 1\\ \hline
\end{array}
\]
}

\subsection{Solutions for $n=29$}

{\scriptsize
\[
\begin{array}{|c|c|c|}
\hline
 m&x& y \\ \hline 
                2& 1& 1\\ \hline
\end{array}
\]
}


\begin{thebibliography}{10}

\normalsize
\baselineskip=17pt

\bibitem{b2010}
A.Bazs\'o, A.B\'erczes, K.Gy\H ory and \'A.Pint\'er,
{\em On the resolution of equations $Ax^n-By^n=C$ in integers $x,y$
and $n\geq 3$}, II, Publ. Math. (Debrecen), {\bf 76}(2010), 227--250.

\bibitem{b2001}
M.A.Bennett, {\em Rational approximation to algebraic numbers of small height: The Diophantine equation $|ax^n -by^n |=1$}, J. Reine Angew. Math. {\bf 535}(2001), 1-49.

\bibitem{maple}
B.W.Char, K.O.Geddes, G.H.Gonnet, M.B.Monagan, S.M.Watt (eds.)
{\em MAPLE, Reference Manual}, Watcom Publications, Waterloo, Canada, 1988.

\bibitem{book}
I.Ga\'al,
{\em Diophantine equations and power integral bases},
Boston, Birkh\"auser, 2002.

\bibitem{gr}
I.Ga\'al and L.Remete,
{\em Binomial Thue equations and power integral bases in pure quartic fields},
JP Journal of Algebra, Number Theory and Applications, {\bf 32}(2014), 49--61.

\bibitem{b2012}
K.Gy\H ory and \'A.Pint\'er, {\em Binomial Thue equations, ternary equations,
and power values of polynomials}, J. Math. Sciences {\bf 180}(2012), 569--580.

\bibitem{pet}
A.Peth\H o, {\em On the resolution of Thue inequalities},
J.Symbolic Comput., {\bf 4}(1987), 103--109.
\end{thebibliography}
\end{document}